\documentclass[12pt,reqno]{amsart}
\usepackage{amssymb,amsmath}
\usepackage{color}
\usepackage{cases}

\def\({\left(}
\def\){\right)}

\def\eb{\varepsilon}

\def\R {\mathbb{R}}

\newcommand{\be}{\begin{equation} }
\newcommand{\ee}{\end{equation} }

\def \rw {\rightarrow}

\def \p {\partial}
\def \px {\partial_{x}}

\def \and{\qquad\text{and}\qquad}

\def\Dt{\partial_t}

\def\({\left(}
\def\){\right)}

\def\eb{\varepsilon}

\def\eb{\varepsilon}
\def\al{\alpha}

\def\la{\lambda}

\def \px {\partial_{x}}
\def \pt {\partial_{t}}

\def\R {\mathbb{R}}

\def \ti {\tilde}
\def \p {\partial}

\def \and{\qquad\text{and}\qquad}

\def\Dt{\partial_t}

\newtheorem{proposition}{Proposition}[section]
\newtheorem{theorem}[proposition]{Theorem}

\newtheorem{lemma}[proposition]{Lemma}
\theoremstyle{definition}

\newtheorem{remark}[proposition]{Remark}

\numberwithin{equation}{section}

\def\be{\begin{equation}}
\def\ee{\end{equation}}
\def\bp{\begin{proof}}
\def\ep{\end{proof}}

\def \no#1#2#3 {{\bf #1} (#3), #2.}
\def \eds#1#2#3 {#1, #2, #3.}

\title[Finite-parameter feedback stabilization of original Burgers' equations]
{Finite-parameter feedback stabilization of original Burgers' equations and Burgers' equation with nonlocal nonlinearities
}

\author[]{ S.~Gumus and V.~K.~Kalantarov }

\address{ (S.~Gumus), Department of Graduate School of Sciences and Engineering\\
\newline\indent Ko{\c c} University, Istanbul}
\address{(V.~K.~Kalantarov) Department of Mathematics
\newline\indent Ko{\c c} University, Istanbul\\
\newline\indent Azerbaijan State Oil and Industry University, Baku
}

\begin{document}

\begin{abstract}
We study the problem of  global exponential stabilization of original Burgers' equations and the Burgers' equation with nonlocal nonlinearities by controllers depending on finitely many parameters. It is  shown that  solutions of the controlled equations are steering a concrete solution of the non-controlled system as $t\rw \infty$ with an exponential rate.
\end{abstract}


\maketitle
\section{Introduction}\label{s0}
We study the initial boundary value problems for the original Burgers' equations (OBE)
\begin{numcases}{}
\Dt v-Uv - \nu\p_x^2v+ 2v \p_xv=0, \ \ \ x\in (0,1),\ t>0,\label{orb1}\\
U'(t)-R + \nu U(t) =-\int_0^1 v^2dx,\   t>0,\ \label{orb2}\\
 v\big|_{x=0}=v\big|_{x=1}=0,   \
U(0)=U_0, \ v\big|_{t=0}=v_0, \label{orbi}
\end{numcases}
and the initial boundary value problem for the viscous Burgers' equation with nonlocal nonlinearity (BNN)
which is closely related to the system \eqref{orb1}, \eqref{orb2} via $U'=0$
\begin{numcases}{}
\Dt v- \nu\p_x^2v + 2v \p_xv-Rv+kv \int_0^1v^2dx =h, \ \ \ x\in (0,1),\ t>0, \label{Bn0}\\
 v\big|_{x=0}=v\big|_{x=1}=0,  \ \
v\big|_{t=0}=v_0. \label{Bn1}
\end{numcases}
Here, the constants $R >0 $ and $\nu >0$ are pressure and kinematic viscosity coefficients, respectively. Burgers introduced the system \eqref{orb1}, \eqref{orb2} as a mathematical model of fluid motion under the action of a constant pressure drop in a straight channel between parallel walls, \cite{Bur1}, and later studied the stabilization of the stationary solutions of this system in detail, \cite{Bur2}. There is a number of papers devoted to the mathematical analysis of Burgers\rq{} original model of turbulence. The problems of  existence and uniqueness of solutions, and the stability of constant solution of the initial boundary value problems for the system \eqref{orb1}, \eqref{orb2} are studied in \cite{Dl1}, \cite{Dl2} and the asymptotic stability of the laminar stationary solution of the initial boundary value problem for this system under homogeneous Dirichlet\rq{}s boundary conditions are proved in \cite{HoOl}. In \cite{Da}, it is shown that the backward in time asymptotic behavior of the solutions to the initial boundary value problem for the system \eqref{orb1}, \eqref{orb2} include complete orbits which grow with an exponential rate as $t\rw \infty$ and  orbits that blow up in a finite time. In \cite{Ede} the author proved the existence of an exponential attractor of the semigroup generated by the initial boundary value problem for the system \eqref{orb1}, \eqref{orb2}, gave an estimate for the dimension of the attractor and  showed that  the dimension of the attractor grew in proportion to $\sqrt{R}\nu^{-1}$. This indicates that the system becomes more unstable as this number grows. The existence of inertial manifolds for Burgers' original mathematical model system of turbulence is investigated in  the paper \cite{IsOh}.

It is the article \cite{DrRe}, where the Burgers' equation with nonlocal nonlinearity \eqref{Bn0}, which is a simplification of the system \eqref{orb1}-\eqref{orb2}, first appeared. Later the asymptotic behavior and the global existence of solutions of the equation \eqref{Bn0} in case of $k<0$, namely unstable Burgers' equation with nonlocal term was studied in \cite{DeKwLe}. The authors showed that the solutions blow up in a finite time under certain conditions on the initial data. Our goal is to study  the problem  of stabilization of solutions by feedback controllers involving finite-dimensional controllers of the problem \eqref{orb1}-\eqref{orbi} and the problem \eqref{Bn0}, \eqref{Bn1}.

There is a vast literature on feedback stabilization and long time behavior of solutions of  the viscous Burgers\rq{} equation
 \be\label{1D} \Dt v- \nu \p_x^2v + v \p_xv=h,\ee and some of its modifications. Some of this literature include  the results on stability of stationary states and existence of finite-dimensional attractors for semigroups generated by  initial boundary value problems for this equation, \cite{Ze}, \cite{ByGySh2}, \cite{Te} , \cite{Ze}; the studies regarding the problem of local stabilization for the viscous Burgers equation, \cite{ByGySh1}, \cite{LyMeTi}, \cite{BaTa} and the results concerning  global stabilization of the viscous Burgers\rq{} equation  using nonlinear Neumann and Dirichlet boundary control laws, \cite{BaKr}, \cite{LiKr2}, \cite{KrMaVa1}, \cite{Fur1}, \cite{FuOs}, \cite{Sm}.

 In this paper, we proved  global exponential stabilization of solutions of \eqref{orb1}, \eqref{orb2} and  global exponential stabilization with a prescribed exponential rate of solutions to \eqref{Bn0} by feedback controllers involving finitely many Fourier modes and finitely many volume elements. Our study is mainly inspired by the recent results  on feedback stabilization by finite-dimensional controllers of nonlinear dissipative PDEs obtained in \cite{AzTi}, \cite{BaWa}, \cite{Che}, \cite{KaOz}, \cite{KaTi}, \cite{LuTi}.
\subsection*{Notation and inequalities}
Throughout this paper we will use the following notations: The inner product in $L^2(0,1)$ is denoted by $(\cdot,\cdot)$ and the corresponding norm by $\|\cdot\|$. The eigenfunctions of the Sturm-Liouville operator $-\frac{d^2}{dx^2}$ under the homogeneous Dirichlet\rq{}s or periodic boundary conditions corresponding to eigenvalues $\la_1,\la_2,\cdots,\la_n,\cdots$ are denoted by $w_1,w_2,\cdots, w_n,\cdots$.

In what follows the nonlinear terms appearing in the calculations are estimated by using the following inequalities:
\begin{itemize}
\item {\it Young\rq{}s inequality:} For each $a, b>0$ and  $\epsilon>0$
\be\label{Young}
ab \le \epsilon \frac{a^p}p+ \frac1{\epsilon^{q/p}}\frac{b^q}q, \ \mbox{where} \  p, q>0 \  \mbox{and} \  \frac1p+\frac1q=1.
\ee
\item {\it Poincar\'e inequality}
\be\label{PF}
\|u\|^2\le \la_1^{-1}\| u'\|^2, \ \ \forall u\in  H_0^1,
\ee
and the inequality
\be\label{PFN}
\sum\limits_{k=N+1}^\infty|(u, w_k)|^2\le \la_{N+1}^{-1}\| u'\|^2, \ \ \forall u\in  H_0^1,
\ee
where  $\la_{n}$ is the $n$'th eigenvalue of the operator   $-\frac{d^2}{dx^2} $  under the homogeneous Dirichlet's boundary conditions.
\item {\it The $1D$ Gagliardo-Nirenberg (GN) inequality}
\be\label{GN}
\|u^{(j)}\|_{L^p(0,1)}\le \beta\|u\|^{1-\al}\|u^{(m)}\|^\al, \ \ \forall u\in H^2(0,1)\cap H_0^m(0,1),
\ee
where $p>2, m=1,2, \ \ \frac jm\le \al\le 1, \ \ \al =\left(\frac12+j-\frac1p\right)m^{-1}.$
\item\begin{lemma}[See \cite{AzTi}]\label{lemma:v1}
Assume that $\phi\in H^1(0,1)$. Then the following inequalities hold true
\begin{align}\label{eqLem1}
\|\phi - \sum_{k=1}^N\bar{\phi}_k \chi_{J_k}\| &\leq h\|\px \phi\|,
\end{align}
where $h= \frac{1}{N}$, $J_k = [(k-1)\frac{1}{N}, k\frac{1}{N})$, for $k=1,2,...,N-1$, $J_N = [\frac{N-1}{N}, 1]$, $\chi_{J_k}$ is the characteristic function of the interval $J_k$ and
\begin{align*}
\bar{\phi}_k = \frac{1}{|J_k|}\int_{J_k}\phi(x)dx.
\end{align*}
\end{lemma}
\end{itemize}

\section{Stabilization of original Burgers' equations by controllers depending on Fourier modes}
In this section we show that strong solutions of the  following feedback system
\begin{numcases}{}
\Dt\ti{v}= \ti{U}\ti{v} + \nu\p_x^2\ti{v}-2\ti{v}\p_x\ti{v} - \mu\sum_{k=1}^N(\tilde{v} - v, w_k)w_k, \  \ x\in(0,1),\label{fdob2}\\
\ti{U}'(t) = R - \nu\ti{U} (t)- \| \ti{v}(t)\|^2, \label{fdob1}\\
\ti{v}\big|_{x=0} = \ti{v}\big|_{x=1} = 0, \ \
\ti{U}(0) = \ti{U}_0, \ \ \ \ti{v}\big|_{t=0}= \ti{v}_0(x).\label{fdob3}
\end{numcases}
are approaching as $t\rw \infty$ the strong solution of the problem \eqref{orb1}-\eqref{orbi} with an exponential rate.\\
A strong solution of the problem \eqref{orb1}-\eqref{orb2} is a pair of functions $[v,U]$ such  that \\
$ v\in L^\infty(0,T; H_0^1(0,1))\cap L^2(0,T; H^2(0,1)), \ U\in C^1(0,T)$ and the equation \eqref{orb1} is satisfied in $L^2$ sense. A strong solution of \eqref{fdob1}-\eqref{fdob3} is defined similarly. We would like to note that the proof of existence and uniqueness of the problems can be done by the standard Galerkin method employing the estimates obtained below (see, e.g., \cite{Dl1}).\\
Let $z := \ti{v} - v$ and $W := \tilde{U}- U$. Then $[z, W]$ is a solution of the system
\begin{numcases}{}
\pt z = \ti{U}z + Wv + \nu\p_x^2z - 2(\tilde{v}\p_xz + z\p_xv) - \mu\sum_{k=1}^N(z,w_k)w_k,\label{nfdob2}\\
W'(t) = -\nu W(t) - \|\ti v(t)\|^2 + \| v(t)\|^2,\label{nfdob1}\\
z\big|_{x=0} = z\big|_{x=1}  = 0, \ \
W(0) = \tilde{U}_0 - U_0, \ \  z\big|_{t=0}= \tilde{v}_0- v_0.\label{nfdob4}
\end{numcases}
Multiplying \eqref{nfdob1} with $W$ and \eqref{nfdob2} by $z$ in  $L^2(0,1)$ and then adding the obtained equations, we get
\begin{multline}\label{est2}
\frac{1}{2}\frac{d}{dt}\left[\|z(t)\|^2 + |W(t)|^2 \right] + \nu\|\p_xz(t)\|^2 + \nu|W(t)|^2 = \ti{U}(t)\|z(t)\|^2\\
 - \mu\sum_{k=1}^N(z,w_k)^2 - W(t)(z,\ti{v})+ (z^2(t), \p_x\ti v(t)-2\p_x v(t)).
\end{multline}
We estimate the last two terms on  the right-hand side of \eqref{est2}, by applying Young's inequality \eqref{Young} and the GN inequality  $\|z\|_{L^4}\le \beta \|z\|^{3/4}\|\p_x z\|^{1/4}$:
\be\label{est3}
|W(t)(z,\ti{v})| \leq \frac{\nu}2|W(t)|^2 + \frac1{2\nu}\|\ti v(t)\|^2\|z(t)\|^2,
\ee
\begin{multline}\label{es*}
|(z^2(t), \p_x \ti v(t)-2\p_x v(t))|\le \|z(t)\|^2_{L^4}\left(\|\p_x\ti v(t)\|+2\|\p_xv(t)\|\right)\\
\le \frac{\nu}4\|\p_xz(t)\|^2 +
\frac34\nu^{-\frac13}\beta^{8/3}\left(\|\p_x\ti v(t)\|+2\|\p_x v(t)\|\right)^{\frac43}\|z(t)\|^2.
\end{multline}
Utilizing the estimates \eqref{est3} and \eqref{es*} in \eqref{est2}, we obtain
\begin{multline}\label{est4}
\frac{1}{2}\frac{d}{dt}\left[\|z(t)\|^2 + |W(t)|^2 \right] + \frac{\nu}2|W(t)|^2 + \frac{3\nu}4\|\p_xz(t)\|^2  \leq - \mu\sum_{k=1}^N|(z,w_k)|^2\\
+ \left(|\ti U(t)|+ \frac1{2\nu}\|\tilde{v}(t)\|^2
 +\frac34\nu^{-\frac13}\beta^{8/3}\left(\|\p_x\ti v(t)\|+2\|\p_xv(t)\|\right)^{\frac43}\right)\|z(t)\|^2.
\end{multline}

We need to show that the terms $\|v(t)\|$, $|U(t)|$, $|\tilde{U}(t)|$, $\|\tilde{v}(t)\|$, $\|\p_x \ti v(t)\|$ and $\|\p_x v(t)\|$ are bounded.
Multiplying \eqref{orb2} by $U$ and \eqref{orb1} by $v$ in $L^2(0,1)$ and then adding the obtained relations we get:
\begin{equation*}
\frac{1}{2}\frac{d}{dt}\left[ \|v(t)\|^2 + |U(t)|^2\right]+ \nu|U(t)|^2 +  \nu\|\p_xv(t)\|^2  = RU(t)\le \frac1{2\nu}{R^2}+ \frac{\nu}2{|U(t)|^2}.
\end{equation*}

Thanks to the  Poincar\'{e} inequality, we obtain
\begin{align*}
\frac{d}{dt}[ \|v(t)\|^2 + |U(t)|^2] + d_0[|U(t)|^2 + \|v(t)\|^2] \leq \frac{1}{\nu}R^2,
\end{align*}
where $d_0 = \nu\min\{1, 2\lambda_1\}$.
The last inequality implies that
\begin{align*}
 \|v(t)\|^2 + |U(t)|^2\leq [\|v_0\|^2 +|U_0|^2 ]e^{- d_0t} + \frac{R^2}{\nu d_0}(1-e^{-d_0 t}).
\end{align*}
Hence, there exists a number $T_1 >0$ such that
\begin{align}\label{est8}
 \|v(t)\|^2 + |U(t)|^2\leq M_1 :=\frac{2R^2}{\nu d_0}, \ \  \forall t\geq T_1.
\end{align}

To get the estimate for  $\|\p_xv(t)\|^2$ first we multiply \eqref{orb1} by $-\p_x^2v$ in $L^2(0,1)$
\be\label{dx3}
\frac{1}{2}\frac{d}{dt}\|\p_xv(t)\|^2 - U(t)\|\p_xv(t)\|^2 + \nu\|\p^2_xv(t)\|^2+\int_0^1(\p_x v(t))^3dx=0.
\ee
Employing the GN inequality $\|(\p_xv)\|_{L^3}\le \beta \|v\|^{\frac5{12}}\|\p_x^2v\|^{\frac7{12}}$
and Young's inequality \eqref{Young} with $\eb=\frac{2\nu}7, \ p=\frac87$, we get
\begin{align}\label{est9}
\Big |\int_0^1(\p_xv(t))^3dx \Big|\le \frac {\nu}4 \|\p_x^2v(t)\|^2+\beta^{24}\nu^{-7}7^72^{-10}\|v(t)\|^{10}.
\end{align}
On the other hand
\begin{equation}\label{est9a}
|U(t)|\|\p_xv(t)\|^2\le |U(t)|\|\p_xv(t)\|\|\p^2_xv(t)\|\le \frac{\nu}4 \|\p^2_xv\|^2+ \frac{1}{\nu}|U(t)|^2\|v(t)\|^2.
\end{equation}
Thus from \eqref{dx3} together with the estimates \eqref{est9} and \eqref{est9a} it follows that
\be\label{est10}
\frac{d}{dt}\|\p_xv(t)\|^2 +\la_1\nu\|\p_xv(t)\|^2\le \frac{1}{\nu}|U(t)|^2\|v(t)\|^2+\beta^{24}\nu^{-7}7^72^{-9}\|v(t)\|^{10}.
\ee
Since  $\|v(t)\|^2, \lvert U\rvert^2(t)\le M_1, \ \ \forall t\ge T_1$ (recall \eqref{est8}), then there exists
$T_2>0$ such that
\begin{align}\label{est13}
\|\p _xv(t)\|^2 &\leq M_2, \ \  \forall t\geq T_2.
\end{align}

In order to get bounds for $|\ti U(t)|$ and $\|\ti v(t)\|$, we multiply the equation
$(\ref{fdob1})$ by $\tilde{U}$ and  $(\ref{fdob2})$ by $\tilde{v}$  in $L^2(0,1)$ and then add the resulting relations:
\begin{align*}
\frac{1}{2}\frac{d}{dt}\left[\|\ti v(t)\|^2 +|\ti U(t)|^2 \right]  &+ \nu|\ti U(t)|^2 + \nu\|\px\ti v(t)\|^2
= R\ti U(t)\nonumber\\&- \mu\sum_{k=1}^N|(\ti v, w_k)|^2
+ \mu\sum_{k=1}^N(\ti v, w_k)(v, w_k).
\end{align*}
By using Young's and Poincar\'{e} inequalities and the estimate \eqref{est8}, we obtain the inequality
\begin{align}\label{est14}
\frac{1}{2}\frac{d}{dt}\left[\|\ti v(t)\|^2 +|\ti U(t)|^2 \right] +\frac{\nu}{2}|\ti U(t)|^2 &+ \la_1\nu\|\ti v(t)\|^2 \leq \frac{R^2}{2\nu} + \frac{\mu}{4}\sum_{k=1}^N(v, w_k)^2\nonumber\\
&\leq \frac{R^2}{2\nu} + \frac{\mu}{4} M_1, \ \ \forall t \geq T_1.
\end{align}
Integrating \eqref{est14}, we get
\be
\|\ti v(t)\|^2 +|\ti U(t)|^2  \leq (\|\ti v(T_1)\|^2 + |\ti U(T_1)|^2 )e^{- d_1(t-T_1)}
+ \frac{1}{d_1}\left(\frac{R^2}{\nu} + \frac{\mu}{2}M_1\right),
\ee
where $d_1 = \nu\min\{1, 2\lambda_1\}$.
Hence, there exists $T_3 > T_1$ such that
\begin{align}\label{est15}
 \|\ti v(t)\|^2 + |\ti U(t)|^2 \leq M_3 := \frac{2R^2}{\nu d_1} + \frac{\mu}{d_1}M_1,\ \  \forall t\geq T_3.
\end{align}
Next, we multiply the equation \eqref{fdob2} by $-\p_x^2\ti v$ in $L^2(0,1)$ and after simple manipulations obtain:
\begin{align*}
\frac{1}{2}\frac{d}{dt}\|\px\ti v(t)\|^2 &- \ti U(t)\|\p_x \ti v(t)\|^2+\nu\|\p^2_x\ti v(t)\|^2+\int_0^1(\p_x \ti v(t))^3dx\\
&=-\mu \sum_{k=1}^N\la_k(\ti v, w_k)^2 +\mu \sum_{k=1}^N\la_k(\ti v, w_k)( v, w_k)\le \frac \mu 4\sum_{k=1}^N\la_k( v, w_k)^2.
\end{align*}
Employing the inequalities \eqref{est9}, \eqref{est9a} we obtain the analog of the inequality \eqref{est10} for $\ti v$
\begin{multline*}
\frac{d}{dt}\|\px \ti v(t)\|^2 +\la_1\nu\|\p_x \ti v(t)\|^2\le \frac{1}{\nu}|U(t)|^2\|\ti v(t)\|^2+\beta^{24}\nu^{-7}7^72^{-9}\|\ti v(t)\|^{10}+\frac\mu2\|\p_xv(t)\|^2,
\end{multline*}
hence, by \eqref{est8}, \eqref{est13} and \eqref{est15} we deduce that
\be\label{tiest}
\| \p_x \ti v(t)\|^2 \le M_4, \ \ \forall t\ge T_4\ge T_3.
\ee
By using the estimates \eqref{est13}, \eqref{est15} and \eqref{tiest}, we infer from  \eqref{est4} the inequality
\begin{align*}
\frac{1}{2}\frac{d}{dt}\left[\|z(t)\|^2 + |W(t)|^2\right] &+ \frac{\nu}{2}|W(t)|^2 + \frac{3\nu}{4}\|\p_xz(t)\|^2 \\
&\leq M_5\|z(t)\|^2 - \mu\sum_{k=1}^N|(z,w_k)|^2,\ \ \forall t \geq T_5,
\end{align*}
where $T_5 := \max\{T_1,T_2, T_3, T_4\}$ and
\begin{equation}\label{M5}
M_5 := \sqrt{M_3}+ \frac1{2\nu}M_3^2 +\frac34\nu^{-\frac13}\beta^{8/3}(\sqrt{M_4} +2  \sqrt{M_2} )^{4/3}.
\end{equation}
Assume that $\mu$ and $N$ are large enough such that
\begin{align}\label{asmp1}
M_5 \leq \mu,\quad \text{and}\quad \lambda_{N+1}^{-1}M_5 \leq \frac{\nu}{4}.
\end{align}
By using these assumptions, the inequality $(\ref{PFN})$ and the Poincar\'{e} inequality, we get
\be\label{zW}
\frac{d}{dt}\left[\|z(t)\|^2 + |W(t)|^2 \right] + \nu|W(t)|^2 + \la_1\nu\|z(t)\|^2 \leq 0, \ \ \forall t \geq T_5.
\ee
From the last inequality we get
\begin{align*}
\|z(t)\|^2 + |W(t)|^2&\leq (\|z(T_5)\|^2 + |W(T_5)|^2 )e^{-d_2 (t - T_5)}, \ \ \forall t \geq T_5,
\end{align*}
where $d_2 = \nu\min\{1, \lambda_1\}$. So we have proved the following
\begin{proposition}\label{theorem1}
Suppose that the conditions $(\ref{asmp1})$ are satisfied. Then there exists $t_0>0$ such that for all $t \geq t_0$ the following inequality holds true
\begin{equation*}
\|\ti v(t)- v(t)\|^2 + |\ti U(t) - U(t)|^2  \leq (\|v(t_0)-\ti v(t_0)\|^2 +|U(t_0)-\ti U(t_0)|^2)e^{-d_2 (t - t_0)},
\end{equation*}
where $d_2 = \nu\min\{1, \lambda_1\}$.
\end{proposition}
Multiplication of \eqref{nfdob2} by $-\p_x^2z$ gives
\begin{multline}\label{H1}
\frac 12\frac d{dt}\|\p_xz\|^2+\nu \|\p^2_xz\|^2=W(t)(\p_xv\ti u,\p_xz)+\ti U\|\p_xz\|^2\\+
2(\ti v\p_xz+z\p_xv,\p^2_xz)-\mu \sum_{k=1}^N\lambda_k(z,w_k)^2.
\end{multline}
Employing the Young inequality and the Sobolev inequality $\|\phi\|^2_{L^\infty}\le c_0\|\p_x\phi\|^2$ we get
\begin{multline}\label{Hb}
2|(\ti v\p_xz+z\p_xv,\p^2_xz)|\le 2\|\ti v\|_{L^\infty}\|\p_xz\|\|\p^2_xz\|+
 2\| z\|_{L^\infty}\|\p_xv\|\|\p^2_xz\|\\
 \le\frac{4c^2_0}{\nu}(\|\p_x \ti v\|^2+\|\p_x \ti v\|^2)\|\p_xz\|^2+\frac \nu2 \|\p^2_xz\|^2.
\end{multline}
and
$$
|W(t)(\p_xv,\p_xz)|\le \frac \nu2 |W(t)|^2+\frac1{2\nu}\|\p_x \ti v\|^2\|\p_x \ti z\|^2
$$
Utilizing last two inequalities in \eqref{H1} we obtain
\be\label{H2}
\frac 12\frac d{dt}\|\p_xz\|^2+\frac\nu2 \|\p^2_xz\|^2\le \frac \nu2 |W(t)|^2+\left[|\ti U|^2+\frac1{2\nu}\|\ti \p_xv\|^2+\frac{4c^2_0}{\nu}(\|\p_x \ti v\|^2+\|\p_x v\|^2)\right]\|\p_xz\|^2.
\ee
Finally adding to \eqref{H2} the inequality \eqref{zW}:
\begin{multline}\label{H3}
\frac d{dt}\left[\frac12 \|\p_xz\|^2+\|z\|^2+|W(t)|^2\right]+\frac \nu2 |W(t)|^2+\lambda_1\nu \|z\|^2+\frac\nu2 \|\p^2_xz\|^2
\\
\le\left[\frac{4\nu c^2_0}\nu M_4+(\frac{4c^2_0}\nu+\frac1{2\nu})M_2+\sqrt{M_3}\right]\|\p_xz\|^2- \mu \sum_{k=1}^N\lambda_k(z,w_k)^2, \ \ \forall t\ge T_5
\end{multline}
Due to the inequality \eqref{PFN} we have
$$
\|\p_xz\|^2=\sum_{k=1}^N\lambda_k(z,w_k)^2+ \sum_{k=N+1}^\infty\lambda_k(z,w_k)^2\le
\sum_{k=1}^N\lambda_k(z,w_k)^2+\lambda_{N+1}^{-1}\|\p^2_xz\|^2.
$$
Thus \eqref{H3} implies

\begin{multline}\label{H4}
\frac d{dt}\left[\frac12 \|\p_xz\|^2+\|z\|^2+|W(t)|^2\right]+\frac \nu2 |W(t)|^2+\lambda_1\nu \|z\|^2\\+ \left(\frac \nu2-\lambda_{N+1}^{-1}Q_0\right)\|\p^2_xz\|^2
\le- (\mu -Q_0)\sum_{k=1}^N\lambda_k(z,w_k)^2,
\end{multline}
where
\be\label{Qu0}
Q_0:=\frac{4\nu c^2_0}\nu M_4+(\frac{4c^2_0}\nu+\frac1{2\nu})M_2+\sqrt{M_3}.
\ee
Suppose that
\be\label{Q0}
\mu -Q_0\ge0, \ \ \lambda_{N+1}^{-1}Q_0\le \frac \nu 4.
\ee
Then from \eqref{H4} we obtain
$$
\frac d{dt}\left[\frac12 \|\p_xz\|^2+\|z\|^2+|W(t)|^2\right]+\alpha_0\left[\frac12 \|\p_xz\|^2+\|z\|^2+|W(t)|^2\right]\le0, \ \ \forall t\ge T_5,
$$
which implies that $\forall t\ge T_5$
\be\label{Hes}
\frac12 \|\p_xz(t)\|^2+\|z(t)\|^2+|W(t)|^2\le \left[\frac12\|\p_xz(T_5)\|^2+\|z(T_5)\|^2+|W(T_5)|^2\right]e^{-(t-T_5)}, \ee
where $\alpha_0=\frac{\nu}2\min\{1,\frac{\lambda_1}2\}$.
Thus we proved the following
\begin{theorem} Suppose that the conditions \eqref{Q0} are satisfied with $Q_0$ defined in \eqref{Qu0}. Then
$\|\p_x\ti{v}(t)-\p_x{v}(t)\|^2+|\ti U(t)-U(t)|^2\rw 0$  with an exponential rate as $t\rw\infty.$
\end{theorem}
\section{ Viscous Burgers' equation with nonlocal nonlinearity}
In this section, we study the problem of stabilization of BNN \eqref{Bn0},
where \\
$h\in L^2(\R^+; L^2(0,1))$ is a given source term, $\nu>0$, $k>0$ and $R>0$ are given numbers.\\
Utilizing the inequalities
\begin{equation*}
|(h,v)|\le \frac{\nu}{2}\|\p_xv\|^2+\frac{1}{2\la_1\nu}\|h\|^2, \ \ R\|v\|^2\le k\|v\|^4+\frac 1{4k}R^2,
\end{equation*}
we obtain from the first energy equality
\begin{equation*}
\frac 12 \frac d{dt}\|v(t)\|^2+ \nu\|\p_x v(t)\|^2-R\|v(t)\|^2+ k\|v(t)\|^4=(h,v),
\end{equation*}
the inequality
$$
\frac d{dt}\|v(t)\|^2+ \nu\la_1\| v(t)\|^2\le \frac 1{2k}R^2+\frac1{\la_1\nu}\|h(t)\|^2.
$$
From this inequality we get the estimate
\begin{equation*}
\|v(t)\|^2\le\|v_0\|^2e^{-\nu\la_1t} +\frac {R^2}{2\la_1\nu k}+\frac1{\la_1\nu}\int_0^t\|h(\tau)\|^2 d\tau
\end{equation*}
which implies that
\be\label{Esv1}
\|v(t)\|\le H_1, \ \ \forall t\ge T_1,
\ee
where $H_1$ depends only on $H_0:= \int_0^\infty\|h(t)\|^2dt.$\\
Next, we multiply \eqref{Bn0} by $-\p_x^2v$ in $L^2(0,1)$ and obtain
\begin{equation}\label{Bnv1}
\frac 12\frac d{dt} \|\p_xv\|^2+ \nu \|\p^2_x v\|^2+\int_0^1(\p_x v)^3dx-R\|\p_xv\|^2+ k\|v\|^2\|\p_xv\|^2=-(h,\p_x^2v).
\end{equation}
Here, we again use  Young\rq{}s inequality \eqref{Young} to obtain:
\be\label{Besv1}
R\|\p_xv\|^2\le R\|v\|\|\p_x^2v\|\le \frac{\nu}8 \|\p^2_x v\|^2+\frac{2R^2}{\nu}\|v\|^2, \ \
|(h,\p_x^2v)|\le \frac{\nu}{8} \|\p^2_x v\|^2+ \frac2{\nu}\|h\|^2.
\ee

Employing the GN inequality and Young's inequality \eqref{Young} with $\eb=\frac{2\nu}{7}, \ p=\frac87$ we get
\begin{align}\label{est9b}
\Big |\int_0^1(\p_xv(t))^3dx \Big|\le \frac{\nu}4\|\p_x^2v(t)\|^2+\beta^{24}7^72^{-10}{\nu}^{-7}\|v(t)\|^{10}.
\end{align}
Utilizing the inequality \eqref{est9b}, two inequalities in \eqref{Besv1} and the Poincar\'{e} inequality, we get from \eqref{Bnv1}
$$
\frac d{dt} \|\p_xv(t)\|^2+ \nu\la_1 \|\p_x v(t)\|^2\le \frac{4}{\nu}\|h(t)\|^2+\frac{4R^2}{\nu}\|v(t)\|^2+\beta^{24}7^72^{-9}{\nu}^{-7}\|v(t)\|^{10}.
$$
Due to \eqref{Esv1} from this inequality we get the estimate
\be\label{esvx}
\|\p_xv(t)\|^2\le H_2, \ \ \forall t\ge T_2,
\ee
where $H_2$ depends  only on $H_0$ and $H_1.$\\
Let us obtain estimates for solutions of the controlled system
\be\label{Bnu}
\begin{cases}
\p_t u- \nu\p_x^2u +2u\px u- Ru+  k\|u\|^2u
= -\mu \sum\limits_{k=1}^N(u-v, w_k)w_k+h,\\
u\Big|_{x=0}=u\Big|_{x=1}=0, \ \ u\Big|_{t=0}=u_0.
\end{cases}
\ee
Multiplication of \eqref{Bnu} by $u$ in $L^2(0,1)$ gives:
\begin{equation}\label{Bnu1a}
\frac 12 \frac d{dt}\|u(t)\|^2+ \nu\|\p_x u(t)\|^2-R\|u(t)\|^2+ k\|u(t)\|^4\le \frac\mu4 \sum\limits_{k=1}^N(v, w_k)^2+(h,u).
\end{equation}
By using the inequalities
\begin{align}\label{BnuEq1}
|(h,u)|\le  \|u\|^2+\frac14\|h\|^2\le \frac k2 \|u\|^4+\frac1{2k}+\frac14\|h\|^2, \ \ R\|u\|^2\le \frac k2 \|u\|^4+\frac1{2k} R^2,
\end{align}
we obtain from \eqref{Bnu1a}
\be\label{Bnu2a}
\frac 12 \frac d{dt}\|u(t)\|^2+ \nu\|\p_x u(t)\|^2
\le
\frac\mu4\|v(t)\|^2+\frac14\|h(t)\|^2+\frac1{2k} (1+R^2).
\ee
We deduce from \eqref{Bnu2a} that
\be\label{esu}
\|u(t)\|\le H_3, \ \forall t\ge T_3,
\ee
where $H_3$ depends on $H_0,H_1,R$ and $k$.
Next, we multiply \eqref{Bnu} by $-\p_x^2u$ in $L^2(0,1)$
\begin{multline}\label{Bnu3a}
\frac12\frac d{dt}\|\p_xu(t)\|^2+\nu\|\p_x^2u(t)\|^2+\int_0^1(\p_xu)^3dx-R\|\p_xu(t)\|^2+k\|u(t)\|^2\|\p_xu(t)\|^2\\
=-\mu \sum\limits_{k=1}^N\la_k(u, w_k)^2+\mu \sum\limits_{k=1}^N\la_k(v, w_k)(u,w_k)-(h,\p_x^2u).
\end{multline}
Thanks to the inequalities \eqref{est9b} and \eqref{Besv1} employed for the term $\px^2u$ and the inequality
$$
\mu\big|\sum\limits_{k=1}^N\la_k(v, w_k)(u,w_k)\big|\le \mu \sum\limits_{k=1}^N\la_k(u,w_k)^2+\frac\mu4 \sum\limits_{k=1}^N\la_k(v,w_k)^2,
$$
\eqref{Bnu3a} implies
$$
\frac12\frac d{dt}\|\p_xu\|^2+\frac{\nu}2\|\p_x^2u\|^2\le \beta^{24}7^72^{-10}{\nu}^{-7}\|u\|^{10}+4R^2 \|u\|^2 +4\|h\|^2+\frac\mu4\|\p_xv\|^2.
$$
Integrating the last inequality we obtain
the next bound for solutions of the problem
\be\label{esux}
\|\px u(t)\|^2\le H_4, \ \ \forall t\ge T_4.
\ee
Finally, we consider the system
\be\label{Bnz}
\begin{cases}
\p_t z-\nu\p_x^2z +2z\p_xu+2v\p_x z-Rz+  k\|u\|^2u-k\|v\|^2v= -\mu \sum\limits_{k=1}^N(z,w_k)w_k,\\
z\Big|_{x=0}=z\Big|_{x=1}=0, \ \ z\Big|_{x=0}=u_0-v_0,
\end{cases}
\ee
where $z=u-v.$
Multiplying the equation \eqref{Bnz} by $z$ in $L^2(0,1)$
and taking into account
$
\left(\|u\|^2u-\|v\|^2v,u-v\right)\ge0,
$
we get
\begin{equation}\label{Bnz1}
\frac 12 \frac d{dt} \|z\|^2+\nu\|\p_xz\|^2+2(z^2, \p_xu)-(z^2, \p_x v) -R\|z\|^2=-\mu \sum\limits_{k=1}^N(z, w_k)^2.
\end{equation}
Then we estimate the third and fourth terms on the left hand side of \eqref{Bnz1}:
\begin{eqnarray}
|(z^2, \p_x v) |\le \|z\|_{L^4}^2\|\p_xv\|&\le&  \beta^2\|z\|^{\frac32}\|\p_xz\|^{\frac12}\|\p_xv\|\nonumber\\
\label{zEq1}&\le& \frac{\nu}4\|\p_xz\|^2+\frac34\beta^{\frac43}{\nu}^{-\frac13}\|z\|^2\|\p_xv\|^{\frac43},
\end{eqnarray}
and similarly
\begin{equation}\label{zEq2}
2|(z^2, \p_x u) |\le \frac{\nu}4\|\p_xz\|^2+\frac34(2\beta)^{\frac43}{\nu}^{-\frac13}\|z\|^2\|\p_xu\|^{\frac43}.
\end{equation}
Thus due to \eqref{esvx} and \eqref{esux}  there exists $T_5\ge T_4$ such that
$$
|(z^2, \p_x v) |+2|(z^2, \p_x u) |\le \frac{\nu}2\|\p_xz\|^2+A_0\|z\|^2, \ \ \forall t\ge T_5,
$$
where \be\label{A0}A_0=\frac34(2\beta)^{\frac43}{\nu}^{-\frac13}H_4^{2/3}+\frac34(\beta)^{\frac43}{\nu}^{-\frac13}H_2^{2/3}.\ee Employing the last inequality we get from \eqref{Bnz1} that
\be\label{Bnz3}
\frac d{dt} \|z(t)\|^2+ \nu\|\px z(t)\|^2-2(A_0+R)\|z(t)\|^2=-2\mu \sum\limits_{k=1}^N(z, w_k)^2.
\ee
Next, we multiply \eqref{Bnz3} by $e^{\sigma t}$ with
$\sigma=\xi-\frac{\la_1\nu}2$ and rewrite the obtained relation in the form
\begin{gather*}
\frac d{dt} \left(e^{\sigma t}\|z(t)\|^2\right)+e^{\sigma t}\nu\|\p_xz(t)\|^2-(\sigma +2A_0+2R)e^{\sigma t}\|z(t)\|^2=-2\mu e^{\sigma t}\sum\limits_{k=1}^N(z, w_k)^2.
\end{gather*}
We rewrite the last equality  in the following form
\begin{multline}\label{Bnzb1}
\frac d{dt} \left(e^{\sigma t}\|z(t)\|^2\right)+e^{\sigma t}\nu\|\p_xz(t)\|^2+[2\mu -(\sigma +2A_0+2R)]e^{\sigma t}
\sum\limits_{k=1}^N(z, w_k)^2\\-
(\sigma +2A_0+2R)e^{\sigma t}\sum\limits_{k=N+1}^\infty(z, w_k)^2=0.
\end{multline}
Since $\sum\limits_{k=N+1}^\infty(z, w_k)^2\le \la_{N+1}^{-1}\|\p_x z(t)\|^2$, we can choose $N$ so large that
$$
\frac{\nu}2 >(\sigma +2A_0+2R)\la_{N+1}^{-1} \ \ \mbox{and} \ \  \mu \ge \frac\sigma2 +R+A_0,$$ and deduce from \eqref{Bnzb1} the  inequality
\be\label{exp1}
\frac d{dt} \left(e^{\sigma t}\|z(t)\|^2\right)+\frac{\nu}2e^{\sigma t}\|\p_xz(t)\|^2\le 0, \ \forall t\ge t_0\ge T_5,
\ee
which implies that
$
\|z(t)\|^2\le e^{-\xi (t-t_0)}\|z(t_0)\|^2.
$\\
So the following proposition holds true
\begin{proposition}
Suppose that $\xi>\frac{\la_1\nu}2$ is an arbitrary number, $N$
 and $\mu$ are so large that
 $$
 \frac{\nu}2 >\la_{N+1}^{-1}(\xi -\frac{\la_1\nu}2+2A_0+2R), \ \  \mu >\frac12 \xi -\frac{\la_1\nu}4+R+A_0,
 $$
 where $A_0$ is defined in \eqref{A0}. Then each strong solution of the problem \eqref{Bnu} is approaching the strong solution of the problem \eqref{Bn0},\eqref{Bn1} with an exponential rate $e^{-\xi t}$ in $L^2(0,1)$ sense.
\end{proposition}
Now we multiply the equation \eqref{Bnz} by $-\p^2_xz$ and obtain:
\begin{multline}\label{Bn1a}
\frac12\frac d{dt}\|\p_xz\|^2+\nu \|\p^2_xz\|^2=2(z\p_xu+v\p_xz,\p^2_xz)+R\|\p_xz\|^2\\
-k\|u\|^2\|\p_xz\|^2-k(u+v,z)(\p_xv,\p_xz)-\mu \sum\limits_{k=1}^N\lambda_k(z, w_k)^2
\end{multline}
By using the Poincar\'e inequality we get
$$
|(u+v,z)(\p_xv,\p_xz)|\le \lambda_1^{-\frac12}(\|u\|+\|v\|)\|\p_xv\|\|\p_xz\|^2,
$$
and similar to \eqref{Hb} we have
$$
2|(z\p_xu+v\p_xz,\p^2_xz)|\le \frac\nu2\|\p^2_xz\|^2+\frac{4c_0^2}\nu(\|\p_xu\|^2+\|\p_xv\|^2)\|\p_xz\|^2.
$$
Utilizing the last two inequalities in \eqref{Bn1a} and employing the estimates \eqref{Esv1},\eqref{esvx}, \eqref{esu}, \eqref{esux} we obtain
$$
\frac12\frac d{dt}\|\p_xz\|^2+\frac\nu2 \|\p^2_xz\|^2\le Q_1\|\p_xz\|^2-\mu \sum\limits_{k=1}^N\lambda_k(z, w_k)^2,
\forall t\ge T_5,$$
where \be\label{Qu1}Q_1:=\frac{4c_0^2}{\nu}(H_4+H_3)+R+k\lambda_1^{-\frac12}(H_1+H_3)\sqrt{H_2}.\ee
Multiplying the last inequality by  $2e^{\sigma t}$ with
$\sigma=\xi-\frac{\la_1\nu}2$, and using the same argumentations when obtaining \eqref{exp1} we get
\be\label{exp1a}
\frac d{dt} \left(e^{\sigma t}\|\p_xz(t)\|^2\right)+\frac{\lambda_1\nu}2e^{\sigma t}\|\p_xz(t)\|^2\le 0, \ \forall t\ge t_0\ge T_5,
\ee
So we proved
\begin{theorem}
If the conditions
$$
\mu>Q_1+\frac12\sigma, \ \ \frac\nu2 >\lambda^{-1}_{N+1}(\sigma+2Q_1)
$$
with $Q_1$ is defined in \eqref{Qu1} are satisfied, then
each solution of the problem \eqref{Bnu} is approaching the strong solution of the problem \eqref{Bn0},\eqref{Bn1} with an exponential rate $e^{-\xi t}$ in $H^1(0,1)$ sense.
\end{theorem}

\begin{remark} Let us note that similarly we can study  the following problem with feedback controller based on finitely many volume elements for the Burgers equation  with nonlocal nonlinearity \eqref{Bn0}   :
\be \label{fVol1}
\begin{cases}
  \pt u-\nu \px^2 u +2u\px u -Ru + ku\|u\|^2 = h -\mu\sum\limits_{k=1}^{N}(\bar{u}_k - \bar{v}_k)\chi_{J_k}(x), \\
u\Big|_{x=0}=u\Big|_{x=1}=0, \ \ u\Big|_{x=0}=u_0.
\end{cases}
\ee
where $J_k, \chi_{J_k}$ and $\bar{u}_k, \bar{u}_k, \ \ k=1,\cdots, N$ are defined in Lemma \ref{lemma:v1}.

Now, by setting $z = u-v$, we see that $z$ is a solution of the following problem
\begin{equation}\label{zfVol1}
\begin{cases}
\pt z - \nu \px^2 z + 2z\px u + 2v\px z -Rz + ku\|u\|^2 - k\|v\|^2v =  - \mu\sum\limits_{k=1}^{N}\bar{z}_k\chi_{J_k}(x),\\
z\Big|_{x=0}=z\Big|_{x=1}=0, \ \ z\Big|_{x=0}=u_0-v_0,
\end{cases}
\end{equation}

Multiplying \eqref{zfVol1} by $z$ in $L^2(0,1)$ and using the inequlity $k(u\|u\|^2 - \|v\|^2v, u-v) \geq 0$, we obtain
\begin{multline}\label{zEq1a}
\frac{1}{2}\frac{d}{dt}\|z\|^2 + \nu\|\px z\|^2 + 2(z\px z,z)+2(v\p_xz+z\p_xv,z)  - R\|z(t)\|^2
\\
\le -\mu  (\sum_{k=1}^N \bar{z}_k\chi_{J_k}(\cdot),z).
\end{multline}
Since
$
(z \px z,z)=0, \ \ 2(z\px z,z)+2(v\p_xz+z\p_xv,z)=(z^2,\p_xv)
$
and
$$
-\mu  (\sum_{k=1}^N \bar{z}_k\chi_{J_k}(\cdot),z)=-\mu  (\sum_{k=1}^N \bar{z}_k\chi_{J_k}(\cdot)-z,z)-\mu \|z\|^2
\le \mu \left\|z-\sum_{k=1}^N \bar{z}_k\chi_{J_k}\right\|\|z\|-\mu \|z\|^2,
$$
utilizing the inequality \eqref{eqLem1} we get from \eqref{zEq1a}
we get
\be\label{Esz1}
\frac{d}{dt}\|z(t)\|^2 + \nu\|\px z(t)\|^2  -R\|z(t)\|^2\le |(z^2,\p_xv)| +\mu h \|\px z(t)\|\|z(t)\|-\mu \|z\|^2.
\ee
Utilizing GN inequality and Young's inequality we have
\be\label{z2b}
|(z^2,\p_xz)|\le \|\p_xv\|\|z\|_{L^4}\le \beta \|\p_xv\|\|z\|^{\frac32}\|\p_xz\|^{\frac12}\\
\le \frac\nu4 \|\p_xz\|^2+\frac43\beta^{\frac43}\nu^{-\frac13}\|\p_xv\|^{\frac43}\|z\|^{2},
\ee
\be\label{z2b}
\mu h \|\px z(t)\|\|z(t)\|\le \frac\nu4\|\px z(t)\|^2+\frac1\nu\mu^2h^2\|z(t)\|^2.
\ee
Therefore thanks to the estimate \eqref{esvx} we deduce from \eqref{Esz1}that
\be\label{Esn2}
\frac{1}{2}\frac{d}{dt}\|z\|^2 + \frac\nu2\|\px z\|^2\le(A_1-\mu)\|z\|^2+\frac1\nu\mu^2h^2\|z(t)\|^2, \ \ \forall t\ge T_2,
\ee
where $
A_1:=R+\frac34\beta^{\frac43}\nu^{-\frac13}H_2^{\frac23}.
$
 Next, we  multiply \eqref{Esn2} by $e^{\sigma t}$ with an arbitrary $\sigma>0$ and after using the Poincare\'e inequality, we get

\be\label{Esn3}
\frac{d}{dt}\left(e^{\sigma t}\|z(t)\|^2\right) +\left[\nu\lambda_1- \frac2\nu\mu^2h^2\right]e^{\sigma t}\|z(t)\|^2
 \le (-2\mu +2A_1+\sigma)\|z(t)\|^2.
\ee
If $\nu\lambda_1\ge \frac4\nu\mu^2h^2 \ \ \mbox{and} \ \ \mu \ge A_1+\frac\sigma2
$ then it follows from \eqref{Esn3} that,
\begin{equation*}
  \|z(t)\|^2 \leq e^{-(\sigma+a_0)(t-t_0)}\|z(t_0)\|^2, \ \ \forall t\geq T_2,
\end{equation*}
where $a_0 = \frac12\nu\lambda_1$. Thus, the following theorem holds true
\begin{theorem}
  Assume that $N=\frac{L}h$ and $\mu$ are large enough so that the conditions
  $$
  \mu\ge \sigma+R+\frac34\beta^{\frac43}\nu^{-\frac13}H_2^{\frac23}\ \ \mbox{and} \ \ \nu\lambda_1\ge \nu \mu^2\frac{4L^2}{N^2},
  $$
  are satisfied. Then each strong solution of the problem \eqref{fVol1} is approaching the solution of \eqref{Bn0}-\eqref{Bn1} with an arbitrary exponential rate in $L^2(0,1)$ sense.
\end{theorem}
\end{remark}

\subsection*{Acknowledgments} We would like to thank Alp Eden and Habiba Kalantarova for  valuable comments.

\end{document}